\documentstyle[12pt]{article}
\begin{document}
\newtheorem{theorem}      {Th\'eor\`eme}[section]
\newtheorem{theorem*}     {theorem}
\newtheorem{proposition}  [theorem]{Proposition}
\newtheorem{definition}   [theorem]{Definition}
\newtheorem{e-lemme}        [theorem]{Lemma}
\newtheorem{cor}   [theorem]{Corollaire}
\newtheorem{resultat}     [theorem]{R\'esultat}
\newtheorem{eexercice}    [theorem]{Exercice}
\newtheorem{rrem}    [theorem]{Remarque}
\newtheorem{pprobleme}    [theorem]{Probl\`eme}
\newtheorem{eexemple}     [theorem]{Exemple}
\newcommand{\preuve}      {\paragraph{Preuve}}
\newenvironment{probleme} {\begin{pprobleme}\rm}{\end{pprobleme}}
\newenvironment{remarque} {\begin{rremarque}\rm}{\end{rremarque}}
\newenvironment{exercice} {\begin{eexercice}\rm}{\end{eexercice}}
\newenvironment{exemple}  {\begin{eexemple}\rm}{\end{eexemple}}
%
%
\newtheorem{e-theo}      [theorem]{Theorem}
\newtheorem{theo*}     [theorem]{Theorem}
\newtheorem{e-pro}  [theorem]{Proposition}
\newtheorem{e-def}   [theorem]{Definition}
\newtheorem{e-lem}        [theorem]{Lemma}
\newtheorem{e-cor}   [theorem]{Corollary}
\newtheorem{e-resultat}     [theorem]{Result}
\newtheorem{ex}    [theorem]{Exercise}
\newtheorem{e-rem}    [theorem]{Remark}
\newtheorem{prob}    [theorem]{Problem}
\newtheorem{example}     [theorem]{Example}
\newcommand{\proof}         {\paragraph{Proof~: }}
\newcommand{\hint}          {\paragraph{Hint}}
\newcommand{\heuristicproof}{\paragraph{heuristic proof}}
\newenvironment{e-probleme} {\begin{e-pprobleme}\rm}{\end{e-pprobleme}}
\newenvironment{e-remarque} {\begin{e-rremarque}\rm}{\end{e-rremarque}}
\newenvironment{e-exercice} {\begin{e-eexercice}\rm}{\end{e-eexercice}}
\newenvironment{e-exemple}  {\begin{e-eexemple}\rm}{\end{e-eexemple}}
\newcommand{\1}        {{\bf 1}}
\newcommand{\pp}       {{{\rm I\!\!\! P}}}
\newcommand{\qq}       {{{\rm I\!\!\! Q}}}
\newcommand{\B}        {{{\rm I\! B}}}
\newcommand{\cc}       {{{\rm I\!\!\! C}}}
\newcommand{\N}        {{{\rm I\! N}}}
\newcommand{\R}        {{{\rm I\! R}}}
\newcommand{\D}        {{{\rm I\! D}}}
\newcommand{\Z}       {{{\rm Z\!\! Z}}}
\newcommand{\C}        {{\bf C}}        
\newcommand{\rank}{\hbox{rank}}
\newcommand{\CC}{{\cal C}}
\def \Re {{\rm Re\,}}
\def \Im {{ \rm Im\,}}
\def\Hom{{\rm Hom\,}}
\def\Lip{{\rm Lip}}
%
%
\newcommand{\dontforget}[1]
{{\mbox{}\\\noindent\rule{1cm}{2mm}\hfill don't forget : #1
\hfill\rule{1cm}{2mm}}\typeout{---------- don't forget : #1 ------------}}
\newcommand{\note}[2]
{ \noindent{\sf #1 \hfill \today}

\noindent\mbox{}\hrulefill\mbox{}
\begin{quote}\begin{quote}\sf #2\end{quote}\end{quote}
\noindent\mbox{}\hrulefill\mbox{}
\vspace{1cm}
}
\title{ Proper J-holomorphic discs in Stein domains of dimension 2 }
\author{ Bernard Coupet{*}, Alexandre Sukhov{**} and Alexander Tumanov{***}}
\date{}
\maketitle

{\small
* Universit\'e de Provence, CMI, 39 rue Joliot-Curie,
Marseille, Cedex, Bernard.Coupet@cmi. math-mrs.fr

** Universit\'e des Sciences et Technologies de Lille, Laboratoire
Paul Painlev\'e,
U.F.R. de
Math\'e-matique, 59655 Villeneuve d'Ascq, Cedex, France,
 sukhov@math.univ-lille1.fr

*** University of Illinois, Department of Mathematics
1409 West Green Street, Urbana, IL 61801, USA,
tumanov@illinois.edu
}
\bigskip

Abstract. We prove the existence of global Bishop discs in a
strictly pseudoconvex Stein domain in an almost complex manifold
of complex dimension 2.
\bigskip


\bigskip

\section{Introduction}

The problem of embedding complex discs or general Riemann
surfaces into complex manifolds has been well-known for a long time.
The interest to the case of almost complex manifolds has grown
due to a strong link with symplectic geometry (Gromov \cite{Gr}).
We present the following result.

\begin{e-theo}
\label{MainTheo}
Let $(M,J)$ be an almost complex manifold of complex dimension 2
admitting a strictly plurisubharmonic exhaustion function $\rho$.
Then for every non-critical value  $c$ of $\rho$, every point
$p \in \Omega_c = \{ \rho  < c \}$ and every vector $v \in T_p(M)$
there exists a $J$-holomorphic immersion
$f:\overline\D\longrightarrow\Omega_c$,
where $\D\subset\cc$ is the unit disc,
such that $f(b\D) \subset b\Omega_c$,
$f(0) = p$, and
$df_0 \left (\frac{\partial}{\partial \Re \zeta}
\right ) = \lambda v$ for some $\lambda >0$.
\end{e-theo}

For a domain $M\subset\cc^n$ with the standard complex structure,
the result is due to Forstneri\v c and Globevnik \cite{FoGl};
there are various generalizations including embedding bordered
Riemann surfaces into singular complex spaces
(see \cite{FoDr} and references there).

Recently Biolley \cite{Bio} proved a similar result for an almost
complex manifold $M$ of any dimension $n$, but under the additional
hypothesis that the defining function $\rho$ is subcritical.
The latter means that $\rho$ does not have critical points of the
maximum Morse index $n$. (A plurisubharmonic function can not
have critical points of index higher than $n$.)
We don't impose such a restriction.
Furthemore, Biolley \cite{Bio} does not prescribe the direction
of the disc. Her method is based on the Floer homology and
substantially uses recent work of Viterbo \cite{Vi}
and Hermann \cite{He}.
Our proof is self-contained; we adapt the ideas of Forstneri\v c and
Globevnik \cite{FoGl} to the almost complex case using the methods
of classical complex analysis and PDE.

In most work on the existence of global discs with boundaries
in prescribed totally real manifolds
(\cite{BeGa, El, Fo1, Hi, Kr1} and others) the authors
use the continuity principle. By the implicit function theorem
and the linearized equation they show that any given disc generates
a family of nearby discs. Then the compactness argument allows
for passing to the limit. In contrast, we construct the discs by
solving the almost Cauchy-Riemann equation directly.

Following \cite{FoGl},
we start with a small disc passing through the given point
in given direction and push the boundary of the
disc in the directions complex-tangent to the level sets of
the defining function $\rho$; it results in increasing $\rho$ due to
pseudoconvexity. This plan leads to a problem of attaching
$J$-holomorphic discs to totally real tori in a level set
of $\rho$. The problem is of independent interest and may
occur elsewhere.  It reduces in turn to
the existence theorem for a boundary value problem for a
quasilinear elliptic system of partial differential equations
in the unit disc (Theorem \ref{Riemann}). We prove it by
the classical methods of the Beltrami equations and
quasiconformal mappings (Ahlfors, Bers, Boyarskii, Lavrentiev,
Morrey, Vekua; see \cite{BJS,Ve} and references there).
The result can be viewed as a far reaching generalization
of the Riemann mapping theorem.

Since the almost Cauchy-Riemann equation is nonlinear,
one can only hope to find a solution close to a current
disc $f$. By measuring the closeness in the $L^p$ norm,
we are able in fact to construct a disc sufficiently far from $f$
in the sup-norm. To make sure we are looking for a disc close
to $f$, we adapt the idea of \cite{FoGl} of adding to $f(\zeta)$
a term with a factor of $\zeta^n$ ($\zeta\in\D$) with big $n$.
We develop a nonlinear version of this idea.

The above procedure works well in the absence of critical
points of $\rho$. In order to push the boundary of the disc
through critical level sets, we use a method by
Drinovec Drnov\v sek and Forstneri\v c \cite{FoDr, Fo},
which consists of temporarily switching
to another plurisubharmonic function at each critical
level set. We point out that adapting this method to
the almost complex case is not a major problem
because the difficulties are localized near the critical
points, in which the almost complex structure can be
closely approximated by the standard complex structure.

Although higher dimension gives one more freedom
for constructing $J$-holomorphic discs,
we must admit that our proof of the main result goes
through in dimension 2 only. The reason is that our
main tool (Theorem \ref{Riemann}) needs a special
coordinate system in which coordinate hyperplanes
$z=const$ are $J$-complex, which generally can be
achieved only in dimension 2.
For a domain in $\cc^n$ with the standard complex structure,
the result is obtained in \cite{FoGl} by reduction
to dimension 2 using sections by 2-dimensional complex
hypersurfaces. Such a reduction in not possible for
almost complex structures.

We thank Franc Forstneri\v c and Josip Globevnik for helpful
discussions, in particular, for pointing out at some difficulties
in the problem and for the important references \cite{FoDr,Fo}.

Parts of the work were completed when the third author was visiting
Universit\'e de Provence and
Universit\'e des Sciences et Technologies de Lille
in the spring of 2006.
He thanks these universities for support and hospitality.

\section{Almost complex manifolds}

Let $(M,J)$ be an almost complex manifold. Denote by $\D$  the unit
disc in $\cc$ and by $J_{st}$  the standard complex structure of
$\cc^n$; the value of $n$ is usually clear from the context.
Let $f$ be a smooth map from $\D$ into $M$. Recall that $f$ is
called {\it $J$-holomorphic}  if $df \circ J_{st} = J \circ df$.
We also call such a map $f$ a $J$-{\it holomorphic} disc or a
{\it pseudoholomorphic} disc or just a {\it holomorphic} disc
when a complex structure is fixed. We will often denote by
$\zeta$ the standard complex coordinate on $\cc$.

A fundamental result of the analysis and geometry of almost complex
structures is the Nijehnuis--Woolf theorem which states that
given point $p \in M$ and given tangent vector $v \in T_pM$ there
exists a $J$-holomorphic disc $f: \D \longrightarrow M$ centered
at $p$, that is, $f(0)=p$ and such that
$df(0)(\partial/\partial \Re \zeta) = \lambda v$
for some $\lambda  > 0$. This disc $f$ depends smoothly on the
initial data $(p,v)$ and the structure $J$. A short proof
of this theorem is given in \cite{Si}. This result will be
used several times in the present paper.

It is well known that an almost complex manifold
$(M,J)$ of complex dimension $n$ can be  locally viewed
as the unit ball $\B$ in $\cc^n$ equipped with an almost
complex structure which is a small deformation of $J_{st}$.
More precisely, let  $(M,J)$ be an almost complex manifold
of complex dimension $n$. Then for every $p \in M$,
$\delta_0 > 0$, and $k\geq 0$ there exist a neighborhood
$U$ of $p$ and a smooth coordinate chart  $z: U \longrightarrow \B$
such that $z(p) = 0$, $dz(p) \circ J(p) \circ dz^{-1}(0) = J_{st}$,
and the direct image $z_*(J) := dz \circ J \circ dz^{-1}$ satisfies
the inequality $\vert\vert z_*(J) - J_{st}
\vert\vert_{C^k(\bar {\B})} \leq \delta_0$.
For a proof we point out that  there exists a diffeomorphism
$z$ from a neighborhood $U'$ of $p \in M$ onto $\B$ such that
$z(p) = 0$ and $dz(p) \circ J(p)\circ dz^{-1}(0) = J_{st}$.
For $\delta > 0$ consider the isotropic dilation
$d_{\delta}: t \mapsto \delta^{-1}t$ in $\cc^n$ and the composite
$z_{\delta} = d_{\delta} \circ z$. Then $\lim_{\delta \rightarrow
0} \vert\vert (z_{\delta})_{*}(J) - J_{st} \vert\vert_{C^k(\bar
{\B})} = 0$. Setting $U = z^{-1}_{\delta}(\B)$ for  positive
$\delta$  small enough, we obtain the desired result.
As a consequence we obtain that
for every point $p\in M$ there exists a neighborhood $U$ of $p$
and a diffeomorphism $z:U \rightarrow \B$ with
center   at $p$ (in the sense that $z(p) =0$) such that the function
$|z|^2$ is $J$-plurisubharmonic on $U$ and
$z_*(J) = J_{st} + O(\vert z \vert)$.

Let $u$ be a  function of class $C^2$ on $M$, let $p \in M$
and $v \in T_pM$. {\it The Levi
form} of $u$ at $p$ evaluated on $v$ is defined by
$L^J(u)(p)(v):=-d(J^* du)(v,Jv)(p)$.

The following result is well known (see, for instance, \cite{DiSu}).
\begin{e-pro}
\label{PROP1}
Let $u$ be a real function of class $C^2$ on $M$, let $p \in M$ and
$v \in T_pM$. Then $ L^J(u)(p)(v) = \Delta(u \circ f)(0)$
where $f: r\D \longrightarrow M$  for some $r > 0$ is an arbitrary $J$-holomorphic map such that
 $f(0) = p$ and $
df(0)(\partial / \partial \Re \zeta) = v$, $\zeta\in r\D$.
\end{e-pro}

The Levi form is  invariant with respect to $J$-biholomorphisms.
More precisely, let $u$ be a $C^2$ real  function on $M$,
let $p \in M$ and $v \in T_pM$.
If $\Phi$ is a $(J,J')$-holomorphic diffeomorphism from
$(M,J)$ into $(M',J')$,
then $L^J(u)(p)(v) =L^{J'}(u \circ \Phi^{-1})(\Phi(p))(d\Phi(p)(v))$.

Finally, it follows from Proposition \ref{PROP1} that a
$C^2$  function $u$ is
$J$-plurisubharmonic on $M$ if and only if $ L^J(u)(p)(v) \geq 0$
for all   $p \in M$, $v \in T_pM$.
Thus, similarly to  the case of
the integrable structure one arrives in a natural way to
the following definition: a $\CC^2$ real valued
function $u$ on $M$ is {\it strictly
$J$-plurisubharmonic} on $M$ if  $ L^J(u)(p)(v)$
is positive for every $p \in M$, $v \in T_pM \backslash \{0\}$.

Let $J$ be a smooth almost complex structure on a neighborhood
of the origin in $\cc^n$ and $J(0) = J_{st}$. Denote by
$z = (z_1,...,z_n)$ the standard coordinates in $\cc^n$
(in matrix computations below we view $z$ as a column).
Then a map $z:\D \longrightarrow \cc^n$ is $J$-holomorphic
if and only if it satisfies the following system of partial
differential equations
\begin{eqnarray}
\label{holomorphy}
z_{\overline\zeta} - A(z)\overline{z_\zeta} = 0,
\end{eqnarray}
where $A(z)$ is the
complex $n\times n$ matrix  defined by
\begin{eqnarray}
\label{matrixA}
A(z)v = (J_{st} + J(z))^{-1}(J_{st} - J(z))\overline v
\end{eqnarray}
It is easy to see that right-hand side of (\ref{matrixA})
is $\cc$-linear in $v\in\cc^n$ with respect to the standard
structure $J_{st}$, hence $A(z)$ is well defined.
Since $J(0) = J_{st}$, we have $A(0) = 0$.
Then in a sufficiently small neighborhood $U$ of the origin
the norm $\parallel A \parallel_{L^\infty(U)}$ is also small,
which implies the ellipticity of the system (\ref{holomorphy}).

However, we will
need a more precise choice of coordinates imposing additional
restrictions on the matrix function $A$.
The proof of the following elementary statement can be found,
for instance, in \cite{DiSu}.

\begin{e-lemme}
\label{normalization}
After a suitable polynomial second degree change of local coordinates near the origin
$$z \mapsto z + \sum a_{kj}z_k\overline z_j$$
we can achieve
\begin{eqnarray*}
A(0) = 0, A_{z}(0) = 0
\end{eqnarray*}
In these coordinates the Levi form of a given $C^2$ function $u$ with respect to
$J$ at the origin coincides with its Levi form with respect to $J_{st}$ that
is
$$L^J(u)(0)(v) = L^{J_{st}}(u)(0)(v)$$
for  every vector $v \in T_0\R^{2n}$.
\end{e-lemme}

\section{Integral transforms in the unit disc}
Let  $\Omega$ be a domain in $\cc$.
Let $T_\Omega$ denote the Cauchy-Green transform
\begin{eqnarray}
T_\Omega f(\zeta) = \frac{1}{2\pi i}\int\int_{\Omega}
\frac{f(\tau) d\tau \wedge d\overline \tau}{\tau - \zeta}.
\end{eqnarray}
Let  $R_\Omega$ denote the Ahlfors-Beurling transform
\begin{eqnarray}
R_\Omega f(\zeta) = \frac{1}{2\pi i}\int\int_{\Omega}
\frac{f(\tau) d\tau \wedge d\overline \tau}{(\tau - \zeta)^2},
\end{eqnarray}
where the integral is considered in the sense of the
Cauchy principal value. We omit the index $\Omega$ if
it is clear form the context.
Denote by $B$ the Bergman projection for $\D$.
$$
Bf(\zeta) = \frac{1}{2 \pi i}\int\int_{\D} \frac{f(\tau)
d\tau \wedge d\overline\tau}{(\overline \tau  \zeta - 1)^2}.
$$
We need the following properties of the above operators.
\begin{proposition}
\label{operators}
\begin{itemize}
\item[(i)]
Let $p > 2$ and $\alpha = (p-2)/p$. Then the linear operator
$T:L^p(\D)\longrightarrow C^\alpha(\cc)$
is bounded, in particular,
$T: L^p(\D) \longrightarrow L^\infty(\D)$ is compact.
If $f\in L^p(\D)$, then
$\partial_{\overline\zeta}Tf=f$, $\zeta\in\D$,
as a Sobolev derivative.
\item[(ii)]
Let $m \geq 0$ be integer and let $0<\alpha<1$.
Then the linear operators
$T:C^{m,\alpha}(\overline\D) \longrightarrow C^{m+1,\alpha}(\cc)$
and
$R:C^{m,\alpha}(\overline \D) \longrightarrow
C^{m,\alpha}(\overline \D)$ are bounded.
Furthermore,if $f\in C^{m,\alpha}(\overline\D)$, then
$\partial_{\overline\zeta}Tf=f$ and
$\partial_{\zeta} Tf = Rf$, $\zeta\in\D$,
in the usual sense.
\item[(iii)]
The operator $R_\Omega$ can be uniquely extended
to a bounded linear operator
$R_\Omega: L^p(\Omega) \longrightarrow L^p(\Omega)$
for every $p > 1$. If $f \in L^p(\D)$, $p > 1$
then $\partial_{\zeta} Tf = Rf$ as
a Sobolev derivative.
Moreover, the operator $R_\cc$ is an isometry of $L^2(\cc)$,
therefore $\parallel R_\cc  \parallel_{L^2(\cc)} = 1$.
\item[(iv)]
The Bergman projection
$B: L^p(\D) \longrightarrow A^p(\D)$ is bounded.
Here $A^p(\D)$ denotes the space of all holomorphic functions
in $\D$ of class $L^p(\D)$.
\item[(v)]
The functions $p \mapsto \parallel T \parallel_{L^p(\Omega)}$
and $p \mapsto \parallel R \parallel_{L^p(\Omega)}$
are logarithmically convex and in particular,
continuous for $p > 1$.
\end{itemize}
\end{proposition}
The proofs of the parts (i)--(iii)  are contained in \cite{Ve}.
The part (iv) follows from (iii); see e. g. \cite{DuShu}.
The part (v) follows by the classical interpolation theorem
of M. Riesz--Torin (see e. g. \cite{Zy}).

We introduce modifications of the operators $T$ and $R$
for solving certain boundary value problems in the unit disc $\D$.
For $f \in L^p(\D)$ we define
\begin{eqnarray}
\label{T0}
T_0f(\zeta) = Tf(\zeta) - \overline{Tf({\overline\zeta\,}^{-1})},
\;\;
\zeta \in \D.
\end{eqnarray}
By Proposition \ref{operators} for $p > 2$ and $\alpha = (p-2)/p$,
the linear operator
$T_0:L^p(\D)\longrightarrow C^\alpha(\D)$
is bounded, in particular,
$T_0: L^p(\D) \longrightarrow L^\infty(\D)$ is compact.
Since the function  $Tf$ is holomorphic and bounded in
$\cc \backslash \overline \D$, then the function
$\zeta \mapsto\overline{(Tf)({\overline\zeta\,}^{-1})}$
is holomorphic in $\D$.
Hence
$\partial_{\overline\zeta} T_0f = \partial_{\overline\zeta} Tf = f$.
Furthermore, for $\zeta\in b\D$, we have
$\zeta={\overline\zeta\,}^{-1}$,
therefore by (\ref{T0}), $\Re T_0f(\zeta)=0$.
Hence for $f\in L^p(\D)$, the function
$u = T_0f$ solves the boundary value problem
$$
\left\{
\begin{array}{cccc}
& \partial_{\overline\zeta} u = f, \zeta \in \D,\\
& \Re u \vert_{b\D} = 0
\end{array}
\right.
$$
We further define
$$
R_0f:= \partial_\zeta T_0f.
$$
Since
$\partial_\zeta Tf = Rf$ and $\partial_{\overline\zeta}Tf = f$,
then
\begin{eqnarray}
\label{R0}
R_0f(\zeta)
=\partial_\zeta T_0f(\zeta)
= Rf(\zeta)
-\overline{\partial_{\overline\zeta} Tf({\overline\zeta\,}^{-1})}
= Rf(\zeta)
+\zeta^{-2}\overline{Rf({\overline\zeta\,}^{-1})},
\end{eqnarray}
and we obtain a nice formula
$$
R_0f=Rf+B\overline f,
$$
where $B$ is the Bergman projection.
By  Propositions \ref{operators}(iv) and (v), the operator
$R_0: L^p(\D) \longrightarrow L^p(\D)$ is bounded, and the map
$p \mapsto \parallel R_0 \parallel_{L^p(\D)}$
is continuous for $p > 1$.
By  Proposition \ref{operators}(iii),
$R$ is an isometry of $L^2(\cc)$.
The analogue of this result for the operator $R_0$
may have been used for the first time by Vinogradov \cite{Vin}.
In fact we came across \cite{Vin} after proving the following

\begin{e-theo}
\label{thmR0}
$R_0$ is a $\R$-linear isometry of $L^2(\D)$,
in particular, $\parallel R_0 \parallel_{L^2(\D)} = 1$.
\end{e-theo}

Since we could not find a proof in the literature,
for completeness we include it here.

\proof
For a domain $G \subset \cc$ we use the inner product
\begin{eqnarray*}
(f,g)_G = -\frac{i}{2}\int \int_G f\overline g
d\zeta \wedge d\overline\zeta.
\end{eqnarray*}
We put
$$
\sigma f(\zeta) = {\overline\zeta\,}^{-2}
\overline {f ({\overline\zeta\,}^{-1})},
\;\;\;\;\;
\psi(\zeta) = \overline\zeta^{2}\zeta^{-2}.
$$
Then $\sigma^2 = {\rm id}$.
By substitution $\zeta \mapsto {\overline\zeta\,}^{-1}$
we obtain
\begin{eqnarray}
\label{inversion}
(\sigma f, \sigma g)_{\D} = (g,f)_{\cc \backslash \D},
\;\;\;\;\;
R \sigma = \psi \sigma R,
\;\;\;\;\;
R  = \psi \sigma R \sigma.
\end{eqnarray}
By (\ref{R0}) we have
$$
R_0f = Rf + \psi \sigma R f.
$$
Let $f \in L^2(\D)$.
Extend $f$ to all of $\cc$ by putting $f(\zeta) = 0$
for $\vert \zeta \vert > 1$.
Then
\begin{eqnarray*}
& &\parallel R_0 f \parallel^2_{L^2(\D)}
= (Rf + \psi \sigma R f, Rf + \psi \sigma R f)_{\D} =\\
& &(Rf,Rf)_{\D} + 2\Re(Rf, \psi \sigma R f)_{\D}
+ (\psi \sigma R f, \psi \sigma R f)_{\D}.
\end{eqnarray*}
Since $|\psi|=1$, by (\ref{inversion})
we obtain
$$
(\psi \sigma R f, \psi \sigma R f)_{\D}
= (\sigma R f,\sigma R f)_{\D}
= (R f, Rf)_{\cc \backslash \D},
$$
$$
(Rf, \psi\sigma Rf)_{\D}
=(\psi\sigma R\sigma f, \psi\sigma Rf)_{\D}
=(R\sigma f, Rf)_{\cc\backslash\D}
=(\psi\sigma Rf, Rf)_{\cc\backslash\D}
=\overline{(Rf, \psi\sigma R f)}_{\cc\backslash\D}.
$$
Then by the previous line and because $R$ is an isometry
\begin{eqnarray*}
2\Re(Rf,\psi \sigma Rf)_{\D}
= \Re(Rf,\psi\sigma Rf)_{\cc}
= \Re(Rf,R\sigma f)_{\cc}
= \Re(f,\sigma f)_{\cc}
= 0.
\end{eqnarray*}
Hence
\begin{eqnarray*}
\parallel R_0 f \parallel^2_{L^2(\D)}
= (Rf,Rf)_{\D} + (Rf, Rf)_{\cc \backslash \D}
=\parallel Rf \parallel^2_{L^2(\cc)}
= \parallel f \parallel^2_{L^2(\cc)}
= \parallel f \parallel^2_{L^2(\D)},
\end{eqnarray*}
which proves the theorem.

\section{ Riemann mapping theorem for an elliptic system}
The Riemann mapping theorem asserts that for every simply connected
domain $G\subset\cc$ there exists a conformal map of $G$ onto $\D$.
If $G$ is smooth, then there is a diffeomorphism
$f:\overline G\longrightarrow\overline\D$, which defines
an almost complex structure $J=f_*(J_{st})$ in $\D$.
Then the Riemann mapping theorem reduces to constructing
a $J$-holomorphic map $z: (\D,J_{st}) \longrightarrow (\D,J)$.
The latter satisfies the Beltrami type equation
$\partial_{\overline\zeta}z = A(z)
\partial_{\overline\zeta}\overline z$, which is equivalent
to the linear Beltrami equaion
$\partial_{\overline z}\zeta + A(z)\partial_z \zeta=0$.
We consider the following more general system
\begin{eqnarray}
\label{mainsystem}
\left\{
\begin{array}{cccc}
& &\partial_{\overline\zeta}z
= a(z,w)\partial_{\overline\zeta}\overline z,\\
& &\partial_{\overline\zeta}w
= b(z,w) \partial_{\overline\zeta}\overline z,
\end{array}
\right.
\end{eqnarray}
which cannot be reduced to a linear one.
Here $z$, $w$ are unknown functions of $\zeta\in\overline \D$ and
$a$, $b$ are $C^\infty$ coefficients.
By eliminating $\zeta$, the system reduces
to a nonhomogeneous quasilinear Beltrami type equation
$\partial_{\overline z}w + a\partial_z w=b$,
but we prefer to deal with (\ref{mainsystem})
directly.

The following theorem is our main technical tool for constructing
pseudoholomorphic discs with boundaries in a prescribed torus.
For $r>0$ denote $\D_r:= r\D$.
\begin{e-theo}
\label{Riemann}
Let $a, b: \overline \D \times \overline \D_{1 +\gamma}
\longrightarrow \cc$ $(\gamma > 0)$ be smooth functions such that
$$
a(z,0) = b(z,0) = 0 \;\;\;\;\;
{\rm and} \;\;\;\;\;
\vert a(z,w) \vert \leq a_0 < 1.
$$
Then there exists $C > 0$ such that for every
integer $n \geq 1$ the system (\ref{mainsystem}) admits
a smooth solution $(z_n,w_n)$ with the following properties:
\begin{itemize}
\item[(i)] $\vert z_n(\zeta) \vert = \vert w_n(\zeta)
\vert = 1$ for $\vert \zeta \vert = 1$.
\item[(ii)] $z_n: \overline \D \longrightarrow \overline \D$
is a diffeomorphism with $z_n(0) = 0$.
\item[(iii)] $\vert w_n(\zeta)\vert \leq C\vert \zeta \vert^n$,
$\vert w_n(\zeta) \vert < 1 + \gamma$.
\end{itemize}
\end{e-theo}

\proof Shrinking $\gamma > 0$ if necessary, we extend the
functions $a$ and $b$
to all of $\cc^2$ preserving their properties.
We will look for a solution of (\ref{mainsystem})
in the form
$$z = \zeta e^u, w = \zeta^n e^v.$$
Then for the new unknowns $u$ and $v$ we have the following boundary
value problem
\begin{eqnarray}
\label{2.2}
\left\{
\begin{array}{cccc}
& &\partial_{\overline\zeta} u
= A(u,v,\zeta)(1 + \overline \zeta \partial_{\overline
  \zeta}\overline u), \zeta \in \D\\
& &\partial_{\overline\zeta} v
= B(u,v,\zeta)(1 + \overline \zeta \partial_{\overline
  \zeta}\overline u), \zeta \in \D\\
& &\Re u(\zeta) = \Re v(\zeta) = 0, \vert \zeta \vert = 1
\end{array}
\right.
\end{eqnarray}
where
\begin{eqnarray*}
& &A = a\zeta^{-1}e^{\overline u - u},\\
& &B = b\zeta^{-n}e^{\overline u - v}.
\end{eqnarray*}
Put $\partial_{\overline\zeta} u = h$ and choose $u$ in the
form $u = T_0h$. Then $\partial_\zeta u = R_0 h$, which we plug
into (\ref{2.2}). We obtain the following system of
singular integral equations for $u$, $v$ and $h$:
\begin{eqnarray}
\label{2.3}
\left\{
\begin{array}{cccc}
& &h = A(1 + \overline \zeta \overline{R_0 h}),\\
& &u = T_0 h,\\
& &v = T_0(B(1 + \overline\zeta \overline{R_0h}))
\end{array}
\right.
\end{eqnarray}
We denote by $\parallel f \parallel_p$
the $L^p$-norm of $f$ in $\D$. Since the function
$p \mapsto \parallel R_0 \parallel_p$ is continuous in $p$ and
$\parallel R_0 \parallel_2 = 1$ we choose $p > 2$ such that
$$a_0 \parallel R_0 \parallel_p < 1.$$
For given $u,v \in L^\infty(\D)$ the map
$h \mapsto A(1 + \overline\zeta\overline{R_0h})$ is a contraction
in $L^p(\D)$ because
$$
\parallel \overline\zeta A \parallel_\infty \parallel
R_0 \parallel_p < 1.
$$
Hence there exists a unique solution $h = h(u,v)$ of the first
equation of (\ref{2.3}) satisfying
\begin{eqnarray}
\label{2.4}
\parallel h \parallel_p \leq \frac{\parallel A \parallel_p}
{1 - a_0 \parallel R_0 \parallel_p}
\end{eqnarray}
Consider the map
$F: L^\infty(\D) \times L^\infty(\D) \longrightarrow L^\infty(\D)
\times L^\infty(\D)$
defined by
\begin{eqnarray*}
F: (u,v) \mapsto (U,V) = (T_0h,
T_0(B(1 + \overline\zeta\overline{R_0h})))
\end{eqnarray*}
where $h = h(u,v)$ is determined above.
Then $F$ is continuous (even Lipschitz) map.
Let
$$
E = \{ (u,v) \in L^\infty (\D) \times L^\infty (\D):
\parallel u \parallel_\infty \leq u_0,
\parallel v\parallel_\infty \leq v_0 \}
$$
We need the following
\begin{e-lemme}
\label{Einvariant}
There exist $u_0>0$, $v_0 > 0$ such that $E$ is invariant under $F$.
\end{e-lemme}

Assuming the lemma, we prove the existence of the solution
of (\ref{2.3}). Indeed, since
$T_0: L^p(\D) \longrightarrow L^\infty(\D)$
is compact for $p > 2$, then $F:E\longrightarrow E$ is compact.
Since $E$ is a bounded, closed and convex,
then the existence of the solution of (\ref{2.3}) follows by
Schauder's principle.
\bigskip

\noindent
{\bf Proof of Lemma \ref{Einvariant} :}
Since $a(z,0) = b(z,0) = 0$, we have
$$
\vert a(z,w) \vert \leq C_1 \vert w \vert,\;\;\;\;\;
\vert b(z,w) \vert \leq C_1 \vert w \vert.
$$
Here and below we denote by $C_j$ constants independent of $n$.
We have
$$
\vert a \vert = \vert a(\zeta e^u,\zeta^n e^v) \vert \leq
C_1 e^{\parallel v \parallel_\infty}\vert \zeta \vert^n,
$$$$
\parallel A \parallel_p = \parallel a \zeta^{-1}\parallel_p \leq
C_2 \parallel \zeta^{n-1}\parallel_p e^{\parallel v \parallel_\infty}
\leq C_3 e^{\parallel v \parallel_\infty}n^{-1/p}.
$$

By (\ref{2.4}),
$\parallel h \parallel_p \leq
C_4 e^{\parallel v \parallel_\infty}n^{-1/p}$, hence
\begin{eqnarray*}
\parallel U \parallel_\infty \leq
C_5 e^{\parallel v \parallel_\infty}n^{-1/p}.
\end{eqnarray*}
Similarly
\begin{eqnarray*}
\vert B \vert = \vert b(\zeta e^u,\zeta^n e^v)\zeta^{-n}
e^{\overline u - v}\vert \leq C_1 e^{\parallel u \parallel_\infty},
\end{eqnarray*}
\begin{eqnarray*}
\parallel B \parallel_\infty \leq C_1 e^{\parallel u \parallel_\infty}
\end{eqnarray*}
\begin{eqnarray*}
\parallel V \parallel_\infty \leq
C_7 (\parallel  B \parallel_p
+ \parallel B \parallel_\infty\parallel h \parallel_p) \leq
C_8 e^{\parallel u \parallel_\infty}.
\end{eqnarray*}
Let $\delta = n^{-1/p}$. Then
\begin{eqnarray*}
& &\parallel U \parallel_\infty \leq C_9 \delta e^{\parallel v \parallel_\infty},\\
& &\parallel V \parallel_\infty \leq C_9 e^{\parallel u \parallel_\infty}
\end{eqnarray*}
Consider the system
\begin{eqnarray*}
u_0 = C_9\delta e^{v_0}, v_0 = C_9 e^{u_0}
\end{eqnarray*}
with the unknowns $u_0$, $v_0$. Then
$$u_0 = C_9 \delta e^{C_9e^{u_0}}$$
For small $\delta > 0$ this equation has two positive roots.
Let $u_0 = u_0(\delta)$ be the smaller root and
$v_0 = v_0(\delta) = C_9e^{u_0}$.
Now if $\parallel u \parallel_\infty \leq u_0$,
$\parallel v \parallel_\infty \leq v_0$, then
$$
\parallel U \parallel_\infty \leq
C_9 \delta e^{\parallel v \parallel_\infty}
\leq C_9 \delta e^{v_0} \leq u_0,
$$
$$
\parallel V \parallel_\infty \leq
C_9 \delta e^{\parallel u \parallel_\infty} \leq
C_9 \delta e^{u_0} \leq v_0
$$
Hence $E$ is invariant under $F$, which proves the lemma.
\bigskip

Thus the solution of (\ref{2.3}) in $L^\infty(\D)$ exists for
$n$ big enough.
Since $h \in L^p(\D)$, $ p > 2$, the second and the third
equations of (\ref{2.3}) imply that
$u,v \in C^\alpha(\overline\D)$, $\alpha=(p-2)/p$.
Since $\partial_{\overline\zeta} u = h \in L^p(\D)$ and
$\partial_\zeta u = R_0 h \in L^p(\D)$ as Sobolev's derivatives,
then $u$ and $v$ are solutions of (\ref{2.2}), hence
$z = \zeta e^u$ and $w = \zeta^n e^v$ are solutions
of (\ref{mainsystem}). By the ellipticity of the system,
$z,w\in C^\infty(\D)$.
The smoothness up to the boundary can be derived directly
from the properties of the Beltrami equation;
it also follows by the reflection principle
for pseudoholomorphic discs attached to totally real
manifolds (see, e.g., \cite{MDS}).

Since the winding number of $z|_{b\D}$ about $0$
equals $1$ and
$\left\vert \partial_{\overline\zeta}z/\partial_\zeta z \right\vert
= \vert a \vert \leq a_0 < 1$ then
$z: \overline\D \longrightarrow \overline\D$ is a homeomorphism
by the classical properties of the Beltrami equation \cite{Ve},
and (ii) follows.

Note that $u_0 \longrightarrow 0$, $v_0 \longrightarrow C_9$ as
$n \longrightarrow \infty$.
Since $T_0:L^p(\D) \longrightarrow C^\alpha(\overline\D)$
is bounded, then we have
$$
\parallel v \parallel_{C^\alpha(\overline\D)} \leq C_{10},\;\;\;
\parallel e^v \parallel_{C^\alpha(\overline\D)} \leq C_{11},
$$
and $|w(\zeta)|\le C_{11} |\zeta|^n$.
Furthermore, since $\vert e^v \vert = 1$ on $b\D$, then
$\vert e^{v(\zeta)}\vert\leq 1+ C_{11}(1-\vert\zeta\vert)^{\alpha}$
for $\vert \zeta \vert < 1$. Then
$|w(\zeta)| \le |\zeta|^n (1+ C_{11}(1-\vert\zeta\vert)^{\alpha})$,
hence $\parallel w \parallel_\infty \longrightarrow 1$ as
$n \longrightarrow\infty$.
Hence $\parallel w \parallel_\infty < 1 + \gamma$ for $n$ big
enough, and (iii) follows.
This completes the proof of Theorem \ref{Riemann}.

\section{ Pseudoholomorphic discs attached to  real tori}

This section concerns the geometrization of Theorem \ref{Riemann}.
We apply Theorem \ref{Riemann} in order to obtain  a crucial
technical result on (approximately) attaching pseudoholomorphic
discs to a given real $2$-dimensional torus in $(M,J)$.
We will use this result later for pushing discs across level sets of
the defining function $\rho$ in Theorem \ref{MainTheo}.

The tori and the discs considered in this section are not arbitrary.
We study a special case which will suffice for the proof of
the main result. Given a psedoholomorphic immersed disc $f$, we
associate with $f$ a real 2-dimensional torus $\Lambda$ formed by
the boundary circles of discs $h_\zeta$ centered at the boundary
points $f(\zeta), \zeta\in b\D$.
Thus, our initial data is a pair $(f,\Lambda)$.
Our goal is to construct a pseudoholomorphic disc with the boundary
attached to the torus $\Lambda$. First we find a suitable
neighborhood of the disc $f$ which can be parametrized by the
bidisc in $\cc^2$. We transport the structure $J$ onto this bidisc
and choose the coordinates there such that the equations for
$J$-holomorphic discs take the form used in Theorem \ref{Riemann}.
The theorem will provide a pseudoholomorphic disc approximately
attached to $\Lambda$.

\subsection{Admissible parametrizations by the bidisc
and generated  tori}
Let $f:\D \longrightarrow (M,J)$ be a $J$-holomorphic disc of class
$C^{\infty}(\overline \D)$. Suppose $f$ is an immersion.
Let $\gamma > 0$.
Given $\zeta \in \overline D$ consider a $J$-holomorphic disc
$$
h_\zeta : (1+\gamma)\D \longrightarrow M
$$
satisfying the  condition $h_\zeta(0) = f(\zeta)$
and such that the direction
$dh_\zeta(0)(\frac{\partial}{\partial \Re\tau})$ is not tangent to $f$.
Admitting some abuse of notation, we sometimes
write $h_{f(\zeta)}$ for $h_\zeta$.

This allows to define a $C^\infty$ map
$$
H: \overline\D \times (1 + \gamma) \overline \D \longrightarrow M,
\;\;\;\;\;
H(\zeta,\tau) = h_\zeta(\tau).
$$
Then $H$ has the following properties:
\begin{itemize}
\item[(i)] For every $\zeta \in \D$ the map $h_\zeta:= H(\zeta,\bullet)$ is $J$-holomorphic.
\item[(ii)] For every $\zeta \in \D$ we have $H(\zeta,0) = f(\zeta)$.
\item[(iii)] For every $\zeta \in \D$ the disc $h_\zeta$ is transversal to $f$ at the point $f(\zeta)$.
\end{itemize}
We assume in addition that
\begin{itemize}
\item[(iv)]
$H: \overline\D \times (1 + \gamma) \overline\D \longrightarrow M$
is locally diffeomorphic.
\end{itemize}

Then $\Lambda = H(b\D \times ( 1 + \gamma) b\D)$ is a real
2-dimensional torus immersed into $M$. It is formed by a family
of topological circles $\gamma_\zeta = h_\zeta((1+ \gamma)b\D)$
parametrized by $\zeta \in b\D$. Every such a circle bounds a
$J$-holomorphic disc $h_\zeta: (1+\gamma)\D \longrightarrow M$
centered at $f(\zeta)$. In particular the torus
$\Lambda$ can be continuously deformed to the circle $f(b\D)$.

If the above conditions (i) - (iv) hold we
say that a map $H$ is an {\it admissible parametrization}
of a neighborhood of $f(\overline\D)$  and $\Lambda$ is
{\it the torus generated by} $H$.

\subsection{Ellipticity of admissible parametrizations}
We prove the
following  consequence of Theorem \ref{Riemann}.

\begin{e-theo}
\label{discintor}
Let $f:\overline\D \longrightarrow (M,J)$ be a
$C^{\infty}$ immersion $J$-holomorphic in $\D$.
Suppose that there exists an admissible parametrization $H$ of a
neighborhood of $f(\D)$  and let $\Lambda$ be the generated torus.
Then there exists an immersed $J$-holomorphic disc
$\tilde f$ of class $C^\infty(\overline\D)$ centered at $f(0)$,
tangent to $f$ at $f(0)$ and
satisfying the boundary condition
$\tilde f(b\D) \subset H(b\D \times b\D)$.
\end{e-theo}
We stress that the boundary of $\tilde f$ is attached to the torus
$H(b\D \times b\D)$ and not to $\Lambda$. However since $\gamma > 0$
can be chosen arbitrarily close to $0$, this leads to the following
result sufficient for applications.

\begin{e-cor}
In the hypothesis of the former theorem for any positive integer
$n$ there exists an immersed $J$-holomorphic disc $ f^n$ of class
$C^\infty(\overline\D)$ centered at $f(0)$, tangent to $f$ at $f(0)$
and such that $dist( f^n(b\D),\Lambda) \longrightarrow 0$ as
$n \longrightarrow \infty$.
\end{e-cor}
Here $dist$ denotes any distance compatible with the topology of $M$.

We begin the proof of Theorem \ref{discintor} with the remark that
the discs $h_\zeta$, $ \zeta \in \overline D$, fill a subset $V$ of
$M$ containing  $f(\overline\D)$ which can be viewed as a fiber
space with the base $f(\overline\D)$ and the generic fiber
$h_\zeta((1+\gamma)\D)$. Therefore the defined above map
$$
H: \overline\D \times (1 + \gamma) \overline\D \longrightarrow V
$$
gives a natural parametrization of $V$ by the bidisc
$U_\gamma := \overline\D \times (1 + \gamma) \overline\D$.
Since $H$ is locally diffeomorphic (see (iv) above)
the inverse map $H^{-1}$ is  defined in a neighborhood of every
point of $V$.  This allows to define the almost complex structure
$\tilde J=H^*(J) = dH^{-1} \circ J \circ dH$ on $U_\gamma$.
The structure $\tilde J$ has a special form.
Indeed, in the standard basis of $\R^4$  we have
\begin{eqnarray}
\label{structure}
\tilde J = \left(
\begin{array}{cll}
\tilde J_{11} & & \tilde J_{12}\\
\tilde J_{21} & & \tilde J_{22}
\end{array}
\right)
\end{eqnarray}
where $\tilde J_{kj}$ are real $2 \times 2$ matrices. We recall that
in this basis the standard complex structure $J_{st}^{(2)}$ of $\cc$
has the form

$$  J_{st}^{(2)} = \left(
\begin{array}{cll}
0 & & -1\\
1 & & 0
\end{array}
\right)
$$
It follows by the property (i)
of $H$ that the maps $\tau \mapsto (c,\tau)$ are
$\tilde J$-holomorphic for every fixed $c$. This implies that
$\tilde J_{12} = 0$ and $\tilde J_{22} = J_{st}^{(2)}$. Furthermore,
since the map $\zeta \mapsto (\zeta, 0)$ is $\tilde J$-holomorphic,
we have
$\tilde J_{11}(z,0) = J_{st}^{(2)}$ and $\tilde J_{21}(z,0) = 0$.

Let now $g:\D \longrightarrow U_\gamma$ be a
$\tilde J$-holomorphic map. If we
set $\zeta = \xi + i \eta$, the
Cauchy--Riemann equations have expressing the
$\tilde J$-holomorphicity of $g$
have the form

\begin{eqnarray}
\label{CR0}
\frac{\partial g}{\partial \xi} + \tilde J \frac{\partial g}{\partial
  \eta} = 0
 \end{eqnarray}
  Suppose now that the matrix $J_{st} + J$ is invertible.
  Then the Cauchy--Riemann equations can be rewritten in the form
\begin{eqnarray}
\label{CR1}
g_{\overline\zeta} + A(g)\overline g_{\overline\zeta} = 0
\end{eqnarray}
where $A$ is defined by (\ref{matrixA}).
If we use the notation $g = (z,w)$, then the Cauchy--Riemann equations
(\ref{CR1}) can be written in the form
\begin{eqnarray}
\label{mainsystem1}
\left\{
\begin{array}{cccc}
& &\partial_{\overline\zeta}z
= a(z,w)\partial_{\overline\zeta}\overline z,\\
& &\partial_{\overline\zeta}w
= b(z,w) \partial_{\overline\zeta}\overline z
\end{array}
\right.
\end{eqnarray}
identical to (\ref{mainsystem}).
Furthermore, since $\tilde J(z,0) = J_{st}$, the conditions
$a(z,0) = b(z,0) = 0$ are satisfied.
\begin{proposition}
\label{ellipticity}
We have $\parallel a \parallel_\infty < 1$.
\end{proposition}
\proof The proof consists of two steps. First we study the
matrix $\tilde J + J_{st}$ which determines the matrix $A$ in the
Cauchy--Riemann equations (\ref{CR1}).
\begin{e-lemme}
\label{nondegeneracy}
The matrix $\tilde J(z,w) + J_{st}$ is non-degenerate for
any $(z,w) \in \overline \D \times (1 + \gamma) \overline \D$.
\end{e-lemme}
\proof It suffices to verify the condition
$\det (\tilde J_{11}(z,w) + J_{st}^{(2)}) \neq 0$.
For every fixed $(z,w)$ the matrix $\tilde J_{11}(z,w)$
defines a complex structure on the euclidean space $\R^2$
so there exists a matrix $P = P(z,w)$ such that
\begin{eqnarray}
\label{matrix}
\tilde J_{11}(z,w) = P J_{st}^{(2)} P^{-1}.
\end{eqnarray}
Recall that the manifold  ${\cal J}_2$ of all complex structures
on $\R^2$ can be identified with the quotient $GL(2,\R)/GL(1,\cc)$
and has two connected components: ${\cal J}_2^+$ and ${\cal J}_2^-$.
A structure $\tilde J_{11}$ belongs to ${\cal J}_2^+$
(resp. to ${\cal J}_2^-$) if in the representation (\ref{matrix})
we have $\det P > 0$ (resp. $\det P <0$).
Suppose now that $\det ( P J_{st}^{(2)} P^{-1} + J_{st}^{(2)}) = 0$
or equivalently $\det (P J_{st}^{(2)}  + J_{st}^{(2)}P) = 0$ at some
point $(z,w)$. If we denote by $p_{jk}$ the entries of the matrix $P$,
the last equality means that $\sum_{jk = 1}^2 p_{jk}^2 = 0$ which
together with the non-degeneracy of $P$ implies that $\det P < 0$
so that $\tilde J_{11}(z,w) \in {\cal J}_2^-$.
On the other hand, for the point $(z,0)$ we have $\det P > 0$
since $\tilde J_{11}(z,0) = J_{st}^{(2)}$ so
$\tilde J(z,0) \in {\cal J}_2^+$.  But we can join the points $(z,0)$
and $(z,w)$ by a real segment, so this contradiction proves lemma.

Now we can conclude the proof of Proposition \ref{ellipticity}.
It follows by Lemma \ref{nondegeneracy} that the Cauchy--Riemann
equations (\ref{CR0}) can be written in the form (\ref{mainsystem1})
on $\overline \D \times (1 + \gamma) \overline \D$.
The Cauchy--Riemann equations are elliptic at every point and this
condition is independent of the choice of the coordinates.
The system (\ref{mainsystem1}) is ellipitic at a point
$(z,w)$ if and only if $\vert a(z,w) \vert \ne 1$.
Since $a(z,0) = 0$ we obtain by connectedness that
$\vert a \vert < 1$ on $\overline \D \times (1 + \gamma) \overline \D$,
which concludes the proof.

Now Theorem \ref{discintor} follows by Theorem \ref{Riemann}.

\subsection{Construction of an admissible parametrization with
a prescribed generated torus}
So far we studied a situation where an admissible parametrization
of a neighborhood of an immersed $J$-holomorphic disc was given
and proved the existence of discs with boundaries close to
the generated torus. In the proof of our main result,
we need an admissible parametrization of a neighborhood of a
$J$-holomorphic disc with a given generated torus.

Let $f: \D \longrightarrow M$ be an immersed
$J$-holomorphic disc of class $C^\infty(\overline\D)$.
We extend $f$ smoothly to a neighborhood of $\overline\D$.
Let $U$ be a small neighborhood of $b\D$.
For every point $f(\zeta)$,  $\zeta\in U$, consider a
$J$-holomorphic disc $h_\zeta: 2\D\longrightarrow M$.
Suppose that the map $h_\zeta$ smoothly depends on $\zeta\in U$.
Thus we obtain a smooth map
$$
H: b\D \times \overline \D \longrightarrow M, \;\;\;\;\;
H: (\zeta,\tau) \mapsto h_\zeta(\tau).
$$
Then ${\Lambda}:= H(b\D \times b\D)$ is a real 2-dimensional torus.
In order to construct an admissible parametrization with the generated
torus $\Lambda$ we need  to extend the map $H$ from the cylinder
$b\D \times \overline \D$ to the bidisc
$\overline \D \times \overline \D$.

\begin{definition}
\label{admissible}
We call the  described above torus $\Lambda$ admissible.
We further put
$X_\zeta:=X_{f(\zeta)}= dh_\zeta(0)(\frac{\partial}{\partial\Re \tau})$
for every $\zeta \in U$.
\end{definition}

\begin{e-theo}
\label{discintor2}
Let $f:\D \longrightarrow (M,J)$ be an immersed $J$-holomorphic
disc of class $C^{\infty}(\overline \D)$.
Let $\Lambda$ be an admissible torus. Then there is
a sequence of admissible tori $\Lambda_n$ converging to $\Lambda$
such that for every $n$ there exists an immersed $J$-holomorphic
disc $ f^n$ of class $C^\infty(\overline\D)$ centered at $f(0)$,
tangent to $f$ at $f(0)$ and satisfying the boundary condition
$f^n(b\D) \subset \Lambda_n$.
\end{e-theo}

\proof
Let $\Lambda$ be an admissible torus and let $X$ be the vector field
given by Definition \ref{admissible}. In general it is impossible to
extend $X$ as a non-vanishing vector field transversal to $f(\D)$
at every point. However, for any integer (not necessarily positive)
$n$ we can consider the discs
$h^n_\zeta:\tau \mapsto h_\zeta(\zeta^n\tau)$,
where $\zeta \in b \D$. Their tangent vectors at the points $f(\zeta)$
are equal to $X^n_\zeta := \zeta^n X_\zeta$,
where by multiplying a vector by a complex number $\zeta^n$
we mean applying the operator $(\Re\zeta+(\Im\zeta) J)^n$.
We need the following

\begin{e-lemme}
After a suitable choice of $n$ the vector field $X_\zeta^n$ can be
extended on the disc as a nonvanishing field transversal to $f$ at
every point.
\end{e-lemme}
\proof
First we look for a global parametrization of a neighborhood of
$f(\overline\D)$. Fix an arbitrary vector field $Y$ transversal to
$f(\overline\D)$ at every point. By Nijenhuis - Woolf theorem we
obtain a family of $J$-holomorphic discs $g_z: w \mapsto g_z(w)$,
$z \in \overline \D$ so that $g_z(0) = f(z)$ and $X_{f(\zeta)}$ is
tangent to $g_z$. Then the map $G: (z,w) \mapsto g_z(w)$ is a
local diffeomorphism from a neighborhood of
$\overline\D \times \overline\D$ onto
a neighborhood of $f(\overline\D)$ and $G(z,0) = f(z)$ so we can use
the coordinates $(z,w)$.  We pull back the vector field  $X$
by $G^{-1}$ and consider the vector field
$(G^{-1})_*(X): \zeta \mapsto (G^{-1})_*(X_\zeta)$. Let $m$ be the
winding number of the $w$-component of the vector field
$(G^{-1})_*(X)$ when $\zeta$ runs along the circle $b\D$.
We set $n = -m$. Then the field $(G^{-1})_*(X^n)$ extends on the
disc $(\zeta,0)$ as a smooth vector field $Z$ transversal to this
disc at every point. Then the map $G_*(Z)$ associates to every point
of $\D$ a vector transversal to $f(\D)$ and so defines the desired
extension $\tilde X^n$ of the vector field $X^n$. This proves the
lemma.

Now by the Nijenhuis - Woolf theorem  there exists  a map
$\tilde h_\zeta :\overline\D\longrightarrow M$ which is
$J$-holomorphic on $\D$ such that $\tilde h_\zeta = h_\zeta$
for every $\zeta$ in a neighborhood of  $b\D$ and the vector
$\tilde X^n_{f(\zeta)}$ is tangent to $h_\zeta$ at the origin.
Thus we can extend $H$ to a function defined on
$\overline \D \times \overline \D$ such that the map
$H(\zeta,\bullet)$ is $J$-holomorphic for any
$\zeta \in \overline \D$. This map $H$ is a local diffeomorphism
and so determines an admissible  parametrization of a neighborhood
of $f(\overline\D)$ such that the generated torus coincides
with $\Lambda$. Theorem \ref{discintor2} now follows
by Theorem \ref{discintor}.

\section{Pushing discs through non-critical levels}
In this section we explain how to push a given disc through
non-critical level sets of a strictly plurisubharmonic function.

\begin{proposition}
\label{noncritical}
Suppose that $\rho$ does not have critical values in the closed
interval $[c_1,c_2]$. Let $f: \D \longrightarrow \Omega_{c_1}$
be an immersed $J$-holomorphic disc such that
$f(b\D)\subset b\Omega_{c_1}$.
Then there exists an immersed $J$-holomorphic disc
$\tilde f: \D \longrightarrow \Omega_{c_2}$ such that
$\tilde f(0) = f(0)$,
$d\tilde f(0) = \lambda df(0)$ for some $\lambda > 0$
and $\tilde f(b\D) \subset b\Omega_{c_2}$.
\end{proposition}

For the proof we need some preparations.
 Let $\rho$ be a strictly plurisubharmonic function on an almost
 complex manifold $(M,J)$. For real $c$ consider the domain
$\Omega_c = \{ \rho < c \}$. Suppose that its boundary has no
critical points.
Let $f: \D \longrightarrow \Omega_c$ be a $J$-holomorphic
disc of class $C^\infty(\overline\D)$ and such that
$f(b\D) \subset b\Omega_c$. For every point $p \in f(b\D)$ consider
a $J$-holomorphic disc $h_p: 2\D\longrightarrow M$ touching
$b\Omega_c$ from outside such that
$\rho \circ h_p\vert_{2\D\backslash\{0\}} > c$.
We call the discs $h_p$ the {\it Levi discs}.
The map $h_p$ can be chosen smoothly depending on $p\in f(b\D)$.

An explicit construction of the Levi discs is given in \cite{FoGl}.
In the almost complex case the proof is similar; the only thing
which has to be justified is the existence of discs $h_p$ touching
a strictly pseudoconvex level set from outside. This was recently
proved by Barraud and Mazzilli \cite{BaMa} and
Ivashkovich and Rosay \cite{IvRo}.
In \cite{DiSu} the result is obtained in any dimension.
For reader's convenience we include a simple proof
(see \cite{DiSu}).

\begin{e-lemme}
\label{tangency}
For a point $p \in b\Omega_c$ there
exists a $J$-holomorphic disc $h_p$ such that $h_p(0) = p$ and
$h_p(\D \backslash \{0 \})$ is contained in
$M \backslash \overline \Omega_c$.
\end{e-lemme}
\proof We fix local coordinates $z = (z_1,z_2)$ near $p$ such that
$p = 0$ and $J(0)=J_{st}$. Denote by $e_j, j = 1,2$ the vectors of
the standard basis of $\cc^2$.
By an additional change of coordinates we may achieve that the map
$h: \zeta\mapsto \zeta e_1$ is $J$-holomorphic on $\D$.
We can assume that the Levi form $L_r^J(0,e_1) = 1$ so that
$$
r(z) = 2\Re z_2+2\Re\sum a_{jk} z_j z_k
+\sum\alpha_{jk}z_jz_k + o(\vert z\vert^2)
$$
with
$$\alpha_{11} = \Delta (r \circ h)(0) = 1.$$
Now for every $\delta  >  0$ consider the non-isotropic dilation
$ \Lambda_\delta:(z_1,z_2) \mapsto
(\delta^{-1/2}z_1,\delta^{-1}z_2)$.
The $J$-holomorphicity of the map $h$ implies that the direct
images $J_\delta:= (\Lambda_\delta)_*(J)$ converge to $J_{st}$
as $\delta \longrightarrow 0$ in the $C^k$ norm for every positive
integer $k$ on any compact subset of $\cc^2$. Similarly, the
functions $r_\delta: = \delta^{-1}r\circ \Lambda^{-1}$ converge
to the function $r_0 := 2\Re z_2 + \vert z_1
\vert^2 + 2 \Re \beta z_1^2$ (for some $\beta \in \cc$).

Consider a $J_{st}$-holomorphic disc
$\hat h: \zeta \mapsto \zeta e_1 - \beta\zeta^2 e_2$.
According to the Nijenhuis-Woolf theorem  for every
$\delta \geq 0$ small enough there exists a $J_\delta$-holomorphic
discs $h^\delta$ such that the family $(h^\delta)_{\delta \geq 0}$
depends smoothly on the parameter $\delta$ and  for every
$\delta\geq 0$ we have
$h^\delta(\zeta) = \zeta e_1 + o(|\zeta|)$ and $h^0 = \hat h$.
Since $(r_0 \circ h^0) (\zeta)= |\zeta|^2$, we obtain that
for $\delta  > 0$ small enough  that
$(r_\delta\circ h^\delta)(\zeta)
= A_\delta(\zeta) + o(|\zeta|^2)$, where
$A_\delta$ is a positive definite quadratic form on $\R^2$.
Since the structures $J_\delta$ and $J$ are biholomorphic,
then the lemma follows.
\bigskip

Thus we obtain a smooth map
$$
H: b\D \times \overline \D \longrightarrow M,
\;\;\;\;\;
H: (\zeta,\tau) \mapsto h_{f(\zeta)}(\tau)=:h_\zeta(\tau).
$$
For simplicity we assume here that $H$ is a local
diffeomorphism although the Levi discs $h_\zeta$ can intersect
even for close values of $\zeta$. We prove in a forthcoming
paper that the pullback $H^*(J)$ of $J$ to the bidisc can
be defined even if $H$ is not a local diffeomorphism.
Thus ${\Lambda}:= H(b\D \times b\D)$ is an admissible
torus and $\rho
\vert_{\Lambda} \geq c + \varepsilon$ for some $\varepsilon > 0$.
We stress that
$\varepsilon$ depends only on $\rho$ (more precisely on a constant
separating the norm of the gradient of $\rho$ from zero) and
the $C^2$-norm of $J$.

Now Theorem \ref{discintor2}
implies that there exists a disc $\tilde f$ with the same direction
as $f$ at the center and with the boundary attached to a torus
arbitrarily close to ${\Lambda}$.  Now we cut off the discs $h_\zeta$
by the level set $\{\rho=c+\varepsilon/2 \}$ and obtain a disc
with boundary attached to this level set.
Indeed, we have the following

\begin{e-lemme}
\label{cut-off}
Suppose that $\rho \circ f \vert_{b\D} \geq c_0$ and $c_0$ is a
non-critical value of $\rho$. Then there exists a $J$-holomorphic
disc $\tilde f$ centered at $f(0)$ and tangent to $f$ at the center
with boundary attached to the level set $\{ \rho = c_0 \}$.
\end{e-lemme}
\proof
By the Hopf lemma the disc $f$ intersects the level set
$\{ \rho = c_0\}$ transversally at every point.
Therefore the open set
$\Omega = \{\zeta\in\D: \rho \circ f(\zeta)<c_0\}$
has a smooth boundary. The set $\Omega$ may be disconnected,
but the connected component of $0\in\Omega$ is simply connected by
the maximum principle applied to the function $\rho \circ f$.
Now the lemma follows via reparametrization by
the Riemann mapping theorem.

Then we again consider the Levi discs for this level set etc.
By iterating this argument a finite number of times we obtain
Proposition  \ref{noncritical}.

\section{Pushing discs through a critical level}

In order to push the boundary of the disc $f$
through critical level sets of $\rho$, we use a method of
\cite{Fo,FoDr}, which consists of temporarily switching
to another plurisubharmonic function at each critical
level set. We need a version of the Morse lemma for almost
complex manifolds.

\begin{proposition}
\label{Morse}
Let $(M,J)$ be an almost complex manifold of complex dimension 2.
Let $\rho$
be a strictly plurisubharmonic Morse function on $M$. Then there
exists
another strictly plurisubharmonic Morse function $\tilde\rho$
close to $\rho$
with the same critical points, such that at each critical point
of Morse index
$k$ in local coordinates given by Lemma \ref{normalization} one has
\begin{eqnarray}
\label{normalform}
\tilde \rho(z) = \tilde\rho(0) + \vert z_1 \vert^2
+ \vert z_2 \vert^2 - a_1 \Re z_1^2 - a_2 \Re z_2^2
\end{eqnarray}
where
\begin{itemize}
\item[(i)] $a_1 = a_2 = 0$ if  $ k = 0$,
\item[(ii)] $a_1 = 2$ and $a_2 = 0$ if  $k = 1$,
\item[(iii)] $a_1 = a_2 = 2$ if   $k = 2$.
\end{itemize}
\end{proposition}

{\bf Remark.} This is a weak version of the Morse lemma because
we change the given function $\rho$ instead of reducing it to
a normal form.
\bigskip

The following result must be well known.
For convenience we include a proof.
\begin{e-lemme}
\label{diagonalization}
Let $B$ be a complex symmetric $n\times n$ matrix.
Then there exists a unitary matrix $U$ such that
$U^tBU$ is diagonal with nonnegative elements.
\end{e-lemme}
\proof
Using coordinate-free language, given a hermitian
positive definite form $H$ and a complex symmetric
bilinear form $B$ on a vector space $V$, $\dim_\cc V=n$,
we need $u_1, ..., u_n\in V$ such that
$$
H(u_i,u_j)=\delta_{ij},\;\;\;\;\;
B(u_i,u_j)=c_i\delta_{ij},\;\;\;
c_i\ge0.
$$
If the above holds with just $c_i\in\cc$, then by rotation
$u_i\mapsto\sigma_i u_i$, $|\sigma_i|=1$, we obtain
$c_i\ge0$.
It suffices to find $u_1\in V$, $H(u_1,u_1)=1$, such that
for every $x\in V$,
\begin{eqnarray}
\label{diagonalization1}
H(x,u_1)=0\;\;\;\;
{\rm implies}\;\;\;\;
B(x,u_1)=0.
\end{eqnarray}
Then the rest of $u_i$ in the $H$-orthogonal complement of $u_1$
are found by induction. Given $u\in V$, by duality, there is
a unique vector $L(u)\in V$ such that
for every $x\in V$,
\begin{eqnarray}
\label{diagonalization2}
H(x,L(u))=B(x,u).
\end{eqnarray}
Then $L:V\to V$ is a $\R$-linear
($\cc$-antilinear) transformation. Since $B$ is symmetric,
then by (\ref{diagonalization2}), $L$ is real symmetric (self-adjoint)
with respect to the form $\Re H$. Then the eigenvalues
of $L$ are real and the eigenvectors are in $V$
(generally they are in $V\otimes_\R\cc$).
Let $u_1\in V$ be an eigenvector of $L$, that is
$L(u_1)=\lambda u_1$, for some $\lambda\in\R$.
We normalize $u_1$ so that $H(u_1,u_1)=1$.
Then for $u=u_1$, (\ref{diagonalization2}) implies
(\ref{diagonalization1}), and the lemma follows.
\bigskip

\noindent
{\bf Proof of Proposition \ref{Morse} :}
Let $p$ be a critical point of $\rho$. Introduce a coordinate
system with the origin at $p$ given by Lemma \ref{normalization}.
In these coordinates the function $\rho$ is strictly
plurisubharmonic at the origin with respect to $J_{st}$. Then
$$
\rho(z) = \rho(0) + \sum a_{ij} z_i \overline z_j
+ \Re \sum b_{ij} z_i z_j + O(\vert z \vert^3),
$$
where $a_{ij} = \overline a_{ji}$ and $b_{ij} = b_{ji}$.
By a linear transformation we can reduce to the form
$a_{ij} = \delta_{ij}$.
If we now make a unitary transformation $z\mapsto Uz$
preserving $|z_1|^2+|z_2|^2$, then the matrix
$B = (b_{ij})$ changes to $U^t B U$.
By Lemma \ref{diagonalization} the expression of
$\rho$ reduces to
$$
\rho(z) = \rho(0) + \vert z_1 \vert^2 + \vert z_2 \vert^2 - \Re
(a_1 z_1^2 + a_2 z_2^2) + O(\vert z \vert^3),
$$
where $a_j \geq 0$, $j=1,2$.
The remainder $\varphi = O(\vert z \vert^3)$ can be removed
by changing $\rho$ to $\tilde \rho = \rho - \varphi \lambda$,
where $\lambda(z) = \lambda_0(z/\varepsilon)$ is a smooth
cut-off function with $\lambda_0 \equiv 1$ in a neighborhood
of the origin and $\lambda_0(z) = 0$ for $\vert z \vert \geq 1$,
$\varepsilon > 0$ small enough.

Since $\varphi(z) = O(\vert z \vert^3)$,
then $\vert d(\varphi \lambda) \vert \leq C \vert z \vert^2$,
$\parallel \varphi \lambda \parallel_{C^2(\cc^2)} \leq C\varepsilon$
where $C > 0$ is independent of $\varepsilon$.
Since $\vert d\rho \vert \geq C \vert z \vert$
in a neighborhood of $0$ for some $C > 0$, then for
small $\varepsilon > 0$ the function $\tilde \rho$ has only one
critical point at the origin, is strictly plurisubharmonic and
matches with $\rho$ for $\vert z \vert > \varepsilon$.

The coefficients $a_j$ can be reduced to the standard values
$0$ and $2$ depending on the index $k$ of the
critical point. We need a cut-off function that falls down
from 1 to 0 sufficiently slowly.

\begin{e-lemme}
\label{cut-off-phi}
Given $\delta > 0$ there exists a smooth non-increasing function $\phi$ with a compact support on  $\R_+$ such that
\begin{itemize}
\item[(i)] $\phi = 1$ near the origin.
\item[(ii)] $\vert t\phi'(t) \vert \leq \delta$.
\item[(iii)] $\vert t^2 \phi''(t) \vert \leq \delta$
\end{itemize}
\end{e-lemme}
The lemma follows because
$\int_1^\infty \frac{dt}{t}=\infty$.
\medskip

Let $b_j = 0$ (resp. $2$) if $0 \leq a_j < 1$ (resp. $a_j > 1$).
Let $\lambda(z) = \phi(\vert z \vert/\varepsilon)$,
where $\phi$ is provided by Lemma \ref{cut-off-phi}
for sufficiently small $\delta$.
Then the function
$$
\tilde \rho(z) = \rho(z) + \lambda
[ (a_1 - b_1) \Re z_1^2 + (a_2 - b_2) \Re z_2^2 ]
$$
for sufficiently small $\varepsilon$ has all the desired properties.
Proposition \ref{Morse} is proved.
\bigskip

Thus in what follows we assume  that $\rho$ has the properties given
by Proposition \ref{Morse}. Let $p$ be a critical point of $\rho$
and $\rho(p) = 0$. Without loss of generality assume that the index
$k$ of $p$ is equal to $1$ or $2$ since the disc obtained
by Proposition \ref{noncritical} cannot approach a minimum of $\rho$.
Choose a small neighborhood $U$ of $p$. By (\ref{normalform})
$\rho$ is strictly plurisubharmonic with respect to $J_{st}$.

We apply the construction of Lemma 6.7 of \cite{Fo}.
Consider $c_0  > 0$ small enough such that $0$ is the only critical
value of $\rho$ in the interval $[-c_0, 3c_0]$. We can assume
that $c_0$ is small enough so that the set
$K(c_0):=\{z:\rho(z)\leq 3c_0,\vert x'\vert^2\leq c_0 \}$
is compactly contained in a  neighborhood of the origin
corresponding to $U$. Here we use the notation $x' = x_1$,
$x'' = x_2$ and $\vert x' \vert^2 = x_1^2$ (resp. $x' = (x_1,x_2)$
and $\vert x' \vert^2 =x_1^2 + x_2^2$ ) if  $k = 1$
(resp.   $ k = 2$). We will use similar notations for the
coordinates  $x$, $y$ and the coordinates   $u$, $v$
introduced below.
Let
\begin{eqnarray}
\label{Etotreal}
E = \{ y' =  0, z'' = 0, \vert x' \vert^2 \leq c_0 \}.
\end{eqnarray}
Then $E$ is a totally real submanifold with boundary and
$\dim E=k$.
Consider the isotropic dilations of coordinates
$$
d_{c_0}: z \mapsto w = u + iv = c_0^{-1/2}z.
$$
Set $J_{c_0} = (d_{c_0})_*(J)$. The structures $J_{c_0}$ converge
to $J_{st}$ in any $C^m$ norm on compact subsets of $\cc^2$
as $c_0 \longrightarrow 0$. Consider the function
$\hat \rho(w) := c_0^{-1} \rho(c_0^{1/2}w)$. This function has
no critical values in $[-1,3]$ and its expression in the
coordinates $w = u + iv$ is the same as the expression
(\ref{normalform}) of  $\rho$ that is
$$\hat \rho(w) = 3v_1^2 + v_2^2 - u_1^2 + u_2^2$$
if $k = 1$ and
$$\hat \rho(w) = 3v_1^2 + 3v_2^2 - u_1^2 - u_2^2$$
if $k = 2$. In particular the set $K = d_{c_0}(K(c_0))$ is given
by $\{ w: \hat \rho(w) \leq 3, \vert u' \vert^2 \leq 1 \}$ and
is a fixed compact independent of $c_0$.

It is important that the origin is a
critical point of the function $\rho$ and the local coordinates
and the function $\rho$ are given by Proposition \ref{Morse}.
This allows to use the isotropic dilations in contrast
with Lemma  \ref{tangency}.

Since  the function $\hat \rho$ is strictly plurisubharmonic
with respect to $J_{st}$, we can apply the
construction of   \cite{Fo} (Lemma 6.7 and section 6.4).
We replace the function $\hat \rho$ by a new
function $\varphi$ defined by
$$\varphi(w) = 3v_1^2 + v_2^2 - h(u_1^2) + u_2^2$$
if $k = 1$ and
$$\varphi(w) = 3v_1^2 + 3v_2^2 - h(u_1^2 + u_2^2)$$
if $k = 2$, where $h\ge0$ is a suitable function.
The construction of $h$ depends on the parameter $c_0$ only. In
our ``delated" coordinates $w$ we apply this construction
taking $c_0 =1$. Namely, according to \cite{Fo} there exist
constants $0 < \tau_0 < \tau_1 < 1$ depending on the eigenvalues
of $\hat \rho$ and a function $\varphi$ strictly plurisubharmonic
on $\cc^2$ with respect to $J_{st}$ satisfying the following
properties:

\begin{itemize}
\item[(i)] $\hat \rho \leq \varphi \leq \hat \rho + \tau_1$,
\item[(ii)] $\hat \rho + \tau_0 \leq \varphi$ on the set
$\{ \vert u' \vert^2 \geq \tau_0 \}$
\item[(iii)] $\varphi = \hat \rho + \tau_1$ on
$\{ \vert u' \vert^2 \geq 1 \}$
\end{itemize}
Since $\hat \rho$ is strictly plurisubharmonic with respect
to the structure $J_{c_0}$, the function $\varphi$ also is
strictly $J_{c_0}$-plurisubharmonic on $\{ \vert u' \vert^2 \geq 1 \}$
in view of (iii). On the other hand the structures $J_{c_0}$
converge to $J_{st}$ in any $C^m$ norm on compact subsets
of $\cc^2$ as $c_0 \longrightarrow 0$. Therefore, since $\varphi$
is strictly $J_{st}$-plurisubharmonic, it also is strictly
$J_{c_0}$-plurisubharmonic on $K$ if $c_0$ is small enough.
Thus, $\varphi$ is strictly $J_{c_0}$-plurisubharmonic
on $\{ \hat \rho \leq 3 \}$.

Now consider the function $\tilde \rho (z) = c_0 \varphi(c_0^{-1/2}z)$
and set $t_0 = \tau_0 c_0$.

The function $\tilde \rho$ satisfies the following properties:
\begin{itemize}
\item[(i)] $\tilde \rho$ is strictly
plurisubharmonic (with respect to $J_{st}$)
in a neighborhood $V \subset U$ of $0$ and $\tilde \rho = \rho + t_1$
on the complement of $V$. Here $t_1 >0$ is a constant.
\item[(ii)] $\tilde \rho$ has no critical values on $(0, 3c_0)$
\item[(iii)] There exists  $t_0 \in (0, c_0)$ such that
\begin{eqnarray}
\label{blowup}
\{ \rho \leq -c_0 \} \cup E \subset
\{ \tilde \rho \leq 0 \} \subset \{ \rho
\leq -t_0 \} \cup E,
\end{eqnarray}
where $E$ is defined above by (\ref{Etotreal}).
\item[(iv)] We have
\begin{eqnarray}
\label{shift}
\{ \rho \leq c_0 \}  \subset \{ \tilde \rho \leq 2c_0 \}
\subset \{ \rho < 3c_0\}
\end{eqnarray}
\end{itemize}

By Proposition \ref{noncritical} we construct an immersed
$J$-holomorphic disc $f$ such that
$-t_0 < \rho \circ f\vert_{b\D} < 0$. The boundary of $f$
is contained in a torus $\Lambda$ formed by discs complex tangent
to a level set of $\rho$.
We will perturb the disc $f$ slightly in order to avoid the
intersection of its boundary with $E$.

\begin{proposition}
\label{discperturb}
Let $f:\overline\D\longrightarrow M$ be an immersed
$J$-holomorphic disc in $(M,J)$, where $\dim_{\cc} M=2$.
Let $E$ be a smooth submanifold in $M$.
Then for every $m \ge2$ there  exists a $J$-holomorphic
disc $\tilde f$ arbitrarily close to $f$ in $C^m(\overline\D)$
such that $\tilde f(0) = f(0)$, $d\tilde f(0) = df(0)$,
and $\tilde f|_{b\D}$ is transverse to $E$.
In particular, if $\dim_{\R} E \leq 2$, then
$\tilde f(b\D)\cap E=\emptyset$.
\end{proposition}
\proof
By the implicit function theorem, the restriction $f|_{bD}$ admits
infinitesimal perturbations in all directions. Then the proposition
follows by the proof of Thom's transversality theorem.
\bigskip

We now assume $f(b\D)\cap E=\emptyset$.
In view of the inclusion (\ref{blowup}) we conclude that
$\tilde\rho > 0$ on $f(b\D)$. By Lemma \ref{cut-off} we cut off
the disc $f$ by a level set $\{ \tilde \rho = c \}$ for some
$c > 0$ to assume that now $f(b\D)$ is contained in this level set.
The function $\tilde \rho$ has no critical values in $(0,3c_0)$.
By Proposition \ref{noncritical} applied to the
disc $f$ and the function $\tilde\rho$ there exists
a new disc $\tilde f$ with the boundary contained in
$\{ \tilde\rho > 2c_0 \}$. In view of (\ref{shift}) we have the
inclusion $\{ \tilde\rho > 2c_0 \} \subset \{ \rho > c_0 \}$.
Now the boundary of $\tilde f$ is outside the critical level
$\{ \rho = 0 \}$ as desired, and we switch back to the
original function $\rho$.

\section{Proof of Theorem \ref{MainTheo}}

Since the function $\rho$ is strictly plurisubharmonic, then after a
generic perturbation of $\rho$ which does not change the given level
set, we can assume that $\rho$ is a Morse function. Let $p$ be the
given point in $D$. If $p$ is not a point of minimum of $\rho$,
we proceed as follows. Consider a small $J$-holomorphic
disc $f$ centered at $p$ with the given direction $v$.
Consider a non-critical level set $\rho = c$ such
that $\rho(p)  < c$. Consider a foliation of a neighborhood
of $f$ by a complex one-parameter family of $J$-holomorphic discs
$h_q$, $q \in f(\D)$ such that the boundaries of these discs are
outside the sublevel set $\rho  < c$.  When  $q$ runs over
the circle $f(b\D)$ these boundaries form a torus. Applying
Proposition \ref{discintor} we obtain a new disc $\tilde f$ centered
at $p$ and still in the same direction at $p$ but with
$\rho \circ \tilde f\vert_{b\D} > 0$.

If $p$ is a point of minimum for $\rho$, we drop this
first step and directly have this situation with $\tilde f = f$.
Now the desired results follow by Proposition \ref{noncritical}
combined with the above argument allowing to push boundaries of
discs through critical levels.


{\footnotesize

\begin{thebibliography}{CIT}


\bibitem{BaMa} J.-F. Barraud, E. Mazzilli,
{\it Regular type of real
hypersurfaces in (almost) complex manifolds},
Math. Z. {\bf 248} (2004), 379--405.


\bibitem{BeGa} E. Bedford, B. Gaveau,
{\it Envelopes of holomorphy of certain $2$-spheres in $\cc^2$},
Amer. J. Math. {\bf 105} (1983), 975--1009.





\bibitem{BJS} L. Bers, F. John, M.Schechter
{\it Partial differential equations}
J.Wiley and Sons, 1964.

\bibitem{Bio} A.-L. Biolley,
{\it Floer homology, symplectic and complex
hyperbolicities}, Preprint, ArXiv math.SG/0404551.

\bibitem{Bo} B. Bojarski,
{\it Generalized solutions of a system of differential equations
of first order of elliptic type with discontinuous coefficients},
Math. Sb. {\bf 43} (1957), 451--503.


\bibitem{DiSu} K. Diederich, A. Sukhov,
{\it Plurisubharmonic exhaustion functions and almost complex Stein
structures}, Preprint, ArXiv math.CV/0603417

\bibitem{FoDr} B.~Drinovec Drnov\v sek, F.~Forstneri\v c,
{\it Holomorphic curves in complex spaces},
Duke Math. J. 139 (2007), 203--254.

\bibitem{DuShu} P. Duren, A. Shuster, {\it Bergman spaces},
Math. Surveys and Monographs, 100, AMS, Providence, RI, 2004.

\bibitem{El} Y.~Eliashberg,
{\it Filling by holomorphic discs and its applications},
London Math. Soc. Lecture Notes, {\bf 151} (1990), 45--67.

\bibitem{Fo1} F.~Forstneri\v c,
{\it Polynomial hulls of sets fibered over the
unit circle},
Indiana Univ. Math. J. {\bf 37} (1988), 869--889.


\bibitem{Fo} F.~Forstneri\v c,
{\it Noncritical holomorphic functions on Stein
manifolds}, Acta Math. {\bf 191} (2003), 143--189.


\bibitem{FoGl} F. Forstneri\v c, J. Globevnik,
{\it Discs in pseudoconvex domains},
Comment. Math. Helv. {\bf 67} (1992), 129--145.

\bibitem{Gr} M. Gromov, {\it Pseudo-holomorphic curves in symplectic
manifolds}, Invent. Math. {\bf 82} (1985), 307--347.

\bibitem{He} D. Hermann,
{\it Holomorphic curves and hamiltonian systems in an
open set with restricted contact type boundary},
Duke Math. J. {\bf 103} (2000),

\bibitem{Hi} R. Hind,
{\it Filling by pseudoholomorphic discs with weakly
pseudoconvex boundary conditions},
GAFA {\bf 7} (1997), 462--495.

\bibitem{IvRo} S. Ivashkovich, J.-P. Rosay,
{\it Schwarz-type lemmas for
solutions of $\overline\partial$-inequalities and complete
hyperbolicity of almost complex structures},
Ann. Inst. Fourier {\bf 54} (2004), 2387--2435.

\bibitem{Kr1} N. Kruzhilin,
{\it Two-dimensional spheres on the boundaries of
    pseudoconvex domains in $\cc^2$},
Izv. Akad. Nauk SSSR Ser. Math. {\bf 52} (1988), 16--40.


\bibitem{MDS} D. McDuff, D. Salamon,
{\it $J$-holomorphic curves and symplectic
topology },
AMS Colloquium Publ., 52, AMS, Providence, RI, 2004.

\bibitem{Si} J.-C. Sikorav,
{\it Some properties of holomorphic
curves in almost complex manifolds}, in ``Holomorphic curves in
Symplectic geometry'',
Ed. M.Audin, J.Lafontane, Birkhauser (1994), 165--189.

\bibitem{SuTu} A. Sukhov, A. Tumanov,
{\it Filling hypersurfaces by discs in
almost complex manifolds of dimension 2},
Indiana Univ. Math. J. {\bf 57} (2008), 509--544.


\bibitem{Ve} I. Vekua,
{\it Generalized analytic functions\/},
{\it Fizmatgiz, Moscow} (1959); English translation:
Pergamon Press, London, and Addison-Welsey, Reading, Massachuset
(1962).

\bibitem{Vin} V. S. Vinogradov,
{\it On a boundary value problem for linear first order elliptic
systems of differential equations in the plane (Russian)},
Dokl. Akad. Nauk SSSR, {\bf 118} (1958), 1059--1062.

\bibitem{Vi} C. Viterbo,
{\it Functors and computations in Floer homology with
applications, Part I}, GAFA, 9 (1999),

\bibitem{Zy} A. Zygmund,
{\it Trigonometric series},
Vol. 2. Cambridge University Press, London 1959.


\end{thebibliography}
}
\end{document}